\documentclass{article}[12pt]

\usepackage{amsmath}
\usepackage{amssymb}
\usepackage{amsthm}
\usepackage[pdftex]{hyperref} 

\newtheorem{theorem}{Theorem}
\newtheorem{lemma}[theorem]{Lemma}
\newtheorem{cor}[theorem]{Corollary}

\newtheorem{claim}[theorem]{Claim}

\newtheorem{Prob}[theorem]{Problem}

\title{\bf
A Note on Upper Bounds\\
for Some Generalized Folkman Numbers\footnote{
Supported by the National Natural Science
Foundation (11361008).}}

\author{Xiaodong Xu\\[-0.1ex]
\small Guangxi Academy of Sciences\\[-0.6ex]
\small Nanning 530007, P.R. China\\[-0.6ex]
\small {\tt xxdmaths@sina.com}\\[1.3ex]\and
Meilian Liang\\[-0.1ex]
\small School of Mathematics and Information Science\\[-0.6ex]
\small Guangxi University, Nanning 530004, P.R. China\\[-0.6ex]
\small {\tt gxulml@163.com}\\[1.3ex]\and
Stanis{\l}aw Radziszowski\\[-0.1ex]
\small Department of Computer Science\\[-0.6ex]
\small Rochester Institute of Technology, Rochester, NY 14623\\[-0.6ex]
\small {\tt spr@cs.rit.edu}\\[3.3ex]
}

\date{\today}

\begin{document}
\maketitle
\thispagestyle{empty}
   
\begin{abstract}
We present some new constructive upper bounds
based on product graphs
for generalized vertex Folkman numbers. They
lead to new upper bounds for some
special cases of generalized edge Folkman numbers,
including $F_e(K_3,K_4-e; K_5) \leq 27$ and
$F_e(K_4-e,K_4-e; K_5) \leq 51$. The latter
bound follows from a construction of a
$K_5$-free graph on 51 vertices, for which
every coloring of its edges with two colors
contains a monochromatic $K_4-e$.
\end{abstract}

\bigskip
\noindent
{\bf Keywords:} Folkman number, Ramsey number\\
{\bf AMS classification subjects:} 05C55, 05C35

\section{Folkman numbers} \label{Folkman}

Let $r, s, a_1, \cdots, a_r$ be integers such
that $r \ge 2$,
$s >  \max \{a_1, \cdots, a_r \}$ and
$\min \{a_1, \cdots, a_r \} \ge 2$.
We write $G \rightarrow (a_1 ,\cdots ,a_r)^v$
(resp. $G \rightarrow (a_1 ,\cdots ,a_r)^e$)
if for every
$r$-coloring of $V(G)$ (resp. $E(G)$), there exists
a monochromatic $K_{a_i}$ in $G$ for some color
$i \in \{1, \cdots, r\}$.
The Ramsey number $R(a_1, \cdots, a_r)$
is defined as the smallest integer $n$ such that
$K_n \rightarrow (a_1, \cdots, a_r)^e$.
The sets of vertex and edge Folkman graphs are defined as
$$\mathcal{F}_v(a_1, \cdots, a_r; s)=\{G \;|\; G \rightarrow
(a_1, \cdots, a_r)^v \textrm{ and } K_s \not\subseteq G\},\textrm{ and}$$
$$\mathcal{F}_e(a_1, \cdots, a_r;s)=\{G\;|\; G \rightarrow
(a_1, \cdots, a_r)^e \textrm{ and } K_s \not\subseteq G\},$$
respectively, and the vertex and edge Folkman numbers
are defined as the smallest orders of graphs in these sets,
namely
$$F_v(a_1, \cdots, a_r;s)=\min\{|V(G)|\;|\;G\in
\mathcal{F}_v(a_1, \cdots, a_r;s) \},\textrm{ and}$$
$$F_e(a_1, \cdots, a_r;s)=\min\{|V(G)|\;|\;G\in
\mathcal{F}_e(a_1, \cdots, a_r;s) \}.$$

The generalized vertex and edge Folkman numbers,
$F_v(H_1, \cdots, H_r;H)$ and $F_e(H_1, \cdots, H_r;H)$,
are defined analogously by considering arrowing graphs
$H_i$ while avoiding $H$, instead of arrowing complete
graphs $K_{a_i}$ while avoiding $K_s$.
The edge Folkman number
$F_e(a_1, \cdots, a_r; k)$ can be seen as a generalization
of the classical Ramsey number $R(a_1, \cdots, a_r)$,
since for $k > R(a_1, \cdots, a_r)$ we clearly have
$F_e(a_1, \cdots, a_r; k) = R(a_1, \cdots, a_r).$

\medskip
In 1970, Folkman \cite{Folkman} proved that
for positive integers $k$ and $a_1, \cdots, a_r$,
$F_v (a_1, \cdots, a_r; k)$  and $F_e (a_1, a_2; k)$
exist if and only if $k > \max \{a_1, \cdots$, $a_r\}$.
Folkman's method did not work for
edge colorings for more than two colors.
The existence of $F_e(a_1, \cdots, a_r; k)$
was proved by Ne\v{s}et\v{r}il and R\"{o}dl
in 1976 \cite{NesetrilRodl}.
Folkman numbers have been studied by many
other authors, in particular in
\cite{drr14,RT,Nkolev2008,LiLin2017a,
LinLi2012a,RodlRS,XuShao2010a,XLSW}.
The current authors studied chromatic variations
of Folkman numbers \cite{XLR0}, and some
existence questions for
$F_v(H_1, \cdots, H_r;H)$ and $F_e(H_1, \cdots, H_r;H)$
\cite{XLR1}.

Perhaps the most wanted Folkman number is $F_e(3,3;4)$,
for which the currently best known bounds are 20
\cite{BikovNenov2016a} and 785 (an unpublished
improvement from 786 obtained by Kauffman, Wickus and the third
author, for more information about the upper bound see
\cite{XLR1}). Further improvements of the
bounds on $F_e(3,3;4)$ seem very difficult,
but some insights can be made into similar
cases involving almost complete graphs $K_k-e$.

For vertex-disjoint graphs $G$ and $H$, their {\em join}
$G+H$ has the vertices $V(G) \cup V(H)$ and edges
$E(G) \cup E(H) \cup E(G,H)$, where $E(G,H)$ is the set
of all possible edges between $V(G)$ and $V(H)$.
Let us also denote $K_k-e$ by $J_k$.
In \cite{XLR1}, we proved
the existence of
$F_e(K_{k+1},K_{k+1};J_{k+2})$ and
$F_v(K_k,K_k;J_{k+1})$, for all $k \ge 3$.
In the same paper we discussed the existence
of some generalized Folkman numbers, especially
in the cases of the form $F_e(K_3,K_3;H)$ for some
small graphs $H$. The latter
includes proofs of nonexistence of the numbers
$F_e(K_3,K_3;J_4)$, $F_e(K_3,K_3;K_2+3K_1)$ and
$F_e(K_3,K_3;K_1+P_4)$, and poses some open cases,
like that for $F_e(K_3,K_3;K_1+C_4)$.

In Section 2 we overview some of the prior constructions
and related upper bounds, and we present our new
constructions. They lead to some new concrete upper bounds,
presented in Section 3, for some special cases including
$F_e(K_3,J_4; K_5) \leq 27$ and
$F_e(J_4,J_4; K_5) \leq 51$.

\section{Constructive Upper Bounds}

The multiplicative upper bound inequality for vertex
Folkman numbers stated in Theorem \ref{K-XLSW} below
was proved independently in \cite{Nkolev2008} and \cite{XLSW}.
A related constructive upper bound for $F_v(k,k;k+1)$
was obtained in \cite{XLSW}, which improved
earlier known bounds, however it is still much weaker
than the best known probabilistic upper bound
for these parameters \cite{drr14}.

\medskip
\begin{theorem} \label{K-XLSW}
{\bf \cite{Nkolev2008,XLSW}}
If $\max\{a_1, \cdots , a_r\} \leq a$ and
$\max\{b_1, \cdots , b_r\} \leq b$, then
$$F_v(a_1 b_1, \cdots, a_r b_r; ab+1) \leq
F_v(a_1, \cdots, a_r; a+1) F_v(b_1, \cdots, b_r; b+1).$$
\end{theorem}

For graphs $G$ and $H$, we will use their
{\em lexicographic product graph} $G[H]$ defined
on the set of vertices $V(G) \times V(H)$ with
$\{(u_1,v_1),(u_2,v_2)\} \in E(G[H])$
if and only if
$\{u_1,u_2\} \in E(G)$ or
($u_1=u_2$ and $\{v_1,v_2\} \in E(H)$).

\medskip
\begin{lemma} \label{prodarrow}
For graphs $G, H$ and $H_i$, and integers $a_j \ge 2$,
$1 \le i,j \le r$,
if $G \rightarrow (a_1, \cdots,a_r)^v$
and $H \rightarrow (H_1, \cdots, H_r)^v$, then
$G[H] \rightarrow (K_{a_1}[H_1], \cdots, K_{a_r}[H_r])^v$.
\end{lemma}

\begin{proof}
Let $G$ and $H$ be any graphs as in the assumptions
of the lemma. Let their sets of vertices be
$U=V(G)$ and $V=V(H)$, respectively, and consider
any partition $V(G[H])=\bigcup_{i=1}^{r}X_i$,
i.e. $r$-coloring $C_v$
of the vertices of $G[H]$. We need to show that for
some color $i$, $1 \le i \le r$, the subgraph induced
by $X_i$ contains $K_{a_i}[H_i]$. Note that for each fixed
$u \in U$, the vertices
$V(u)=\{(u,v)\;|\;v \in V\}$ induce a graph
isomorphic to $H$ in $G[H]$. Hence, for each $u \in U$
there exists a color $i(u)$,
$1 \le i(u) \le r$, such that
the subgraph induced by $V(u)$ contains $H_{i(u)}$
in color $i(u)$.

Next, consider the $r$-coloring $C_v^{'}$
of vertices of $G$ defined by $i(u)$.
Since $G \rightarrow (a_1, \cdots,a_r)^v$,
then there exists $j$ such $C_v^{'}$
contains $K_{a_j}$ in color $j$ in $G$, or equivalently,
in the vertex $r$-coloring $C_v$
of $G[H]$ we have $a_j$ isomorphic copies of $H$,
each of them containing $H_j$ in color $j$,
and all of them are interconnected by edges in $G[H]$.
\end{proof}

\medskip
Let $cl(H)$ denote the clique number of graph $H$,
i.e. the largest integer $s$ such that $K_s \subset H$.
The following generalizes Theorem \ref{K-XLSW}.

\medskip
\begin{theorem}  \label{FolkmangeneralGH}
If $\max \{a_1, \cdots, a_r\} \leq a$ and
$\max \{cl(H_1), \cdots, cl(H_r)\} \leq b$, then
$$F_v(K_{a_1}[H_1], \cdots, K_{a_r}[H_r];ab+1)
\leq F_v(a_1, \cdots, a_r; a+1) F_v(H_1, \cdots, H_r;b+1).$$
\end{theorem}

\begin{proof}
Consider any graph
$G \in \mathcal{F}_v(a_1,\cdots,a_r;a+1)$
with the set of vertices $V(G)=U=\{u_1,\cdots,u_s\}$,
where $s=F_v(a_1,\cdots,a_r;a+1)$,
and any graph $H$ such that
$H \in \mathcal{F}_v(H_1,\cdots,H_r;b+1)$
and $V(H)=V=\{v_1,\cdots,v_t\}$,
where $t=F_v(H_1,\cdots,H_r;b+1)$.
Note that $st=|V(G[H])|$ is also equal to
the right hand side of the target inequality.
By the construction of $G[H]$
one can easily see that $cl(G[H]) \le ab$. Finally,
Lemma \ref{prodarrow} implies that
$G[H] \rightarrow (K_{a_1}[H_1], \cdots, K_{a_r}[H_r])^v$,
which completes the proof.
\end{proof}

\medskip
We note now, and will also use it later, that
Theorem \ref{FolkmangeneralGH} is specially
interesting in the cases involving graphs $J_k$,
because we can use the fact that $J_{sk+1}$
is a subgraph of $K_s[J_{k+1}]$. For instance,
using Theorem \ref{FolkmangeneralGH} for two
colors with $s=a_1=a_2=2$, $k=b$, $F_v(2,2;3)=5$
and $H_1=H_2=J_{k+1}$, we obtain
$$F_v(K_{2k+2}-2K_2,K_{2k+2}-2K_2;2k+1)
\leq 5F_v(J_{k+1},J_{k+1};k+1).$$

\smallskip
Further, since
$J_{2k+1}$ is a subgraph of $K_{2k+2}-2e$,
it also holds that
$$F_v(J_{2k+1},J_{2k+1};2k+1) \leq
5F_v(J_{k+1},J_{k+1};k+1).$$

\medskip
In fact, we can do a little better on 3 out of 5
blocks of $F_v(J_{k+1},J_{k+1};k+1)$ vertices,
as stated in the next theorem.

\medskip
\begin{theorem}  \label{FvK{2k+1}-e}
For every integer $k \geq 2$, we have that
$F_v(J_{2k+1},J_{2k+1};2k+1) \leq
2F_v(k,k;k+1) + 2F_v(J_{k+1},J_{k+1};k+1)
+ F_v(K_k,J_{k+1};k+1).$
\end{theorem}

\begin{proof}
Consider any graphs $H_1, H_2, H_3$ such that
$H_1 \in \mathcal{F}_v(k,k;k+1)$,
$H_2 \in \mathcal{F}_v(K_k,J_{k+1};k+1)$, and
$H_3 \in \mathcal{F}_v(J_{k+1},J_{k+1};k+1)$,
and they have the smallest possible number of vertices,
i.e. $|V(H_1)| = F_v(k,k;k+1)$,
$|V(H_2)| = F_v(K_k,J_{k+1};k+1)$, and
$|V(H_3)| = F_v(J_{k+1},J_{k+1};k+1)$.
Let $H_4$ be an isomorphic copy of $H_1$,
and $H_5$ an isomorphic copy of $H_3$.
The clique number of all graphs $H_i$ is
equal to $k$.

Our goal is to construct graph
$G \in \mathcal{F}_v(J_{2k+1},J_{2k+1};2k+1)$
on the set of vertices
$V(G)=\bigcup_{i=1}^5 V(H_i)$,
which has $cl(G)=2k$.
This will suffice to complete the proof.
The set of edges of graph $G$ is defined by
$E(G) = \bigcup_{i=1} ^5 E(H_i) \cup E(1,3)
\cup E(1,5) \cup E(2,4) \cup E(2,5) \cup E(3,4)$,
where
$E(i,j) = \{ \{ u,v\} \;|\; u \in V(H_i), v \in V(H_j),
i \not= j, 1 \le i,j \le 5\}$.
One can easily check that $G$ is $K_{2k+1}$-free,
since the edges of types $E(i,j)$ do not form
any triangles. It remains to be shown that
$G \rightarrow (J_{2k+1},J_{2k+1})^v$.

For a contradiction, suppose that there exists a partition
$V(G)=R \cup B$, i.e. a red-blue coloring of the vertices
of $G$, which has no monochromatic $J_{2k+1}$.
Without loss of generality we may assume that
$H_1$ contains a red $K_k$. Therefore, there is no red
$J_{k+1}$ in $H_3$ and no red $J_{k+1}$ in $H_5$,
otherwise we would have a red $J_{2k+1}$.
Hence, there are blue $J_{k+1}$'s in both
$H_3$ and $H_5$.
Therefore, there is no blue $K_k$ in $H_2$
and no blue $K_k$ in $H_4$.
Hence, there is a red $J_{k+1}$ in $H_2$
and a red $K_k$ in $H_4$, and together
they form a red $J_{2k+1}$.
Thus $G \rightarrow (J_{2k+1},J_{2k+1})^v$,
which completes the proof.
\end{proof}

\smallskip
As an application of the last theorem
we consider an interesting case of $F_v(J_5,J_5;5)$.
Likely, it is just somewhat larger (and harder to compute)
than the well studied classical case of
$F_v(4,4;5)$, for which the currently best
known bounds are
$17 \le F_v(2,3,4;5) \le F_v(4,4;5) \le 23$ \cite{Fv445}.

\begin{claim}
$F_v(J_5,J_5;5) \le 36$.
\end{claim}

\begin{proof}
The first three graphs $H_i$ in the proof
of Theorem \ref{FvK{2k+1}-e}
for this case are of orders equal to
$F_v(2,2;3)$, $F_v(J_3,J_3;3)$ and
$F_v(K_2,J_3;3)$, respectively. The first of them,
equal to 5, is uniquely witnessed by the cycle
$C_5$. Obtaining upper bounds for the other two
requires some work when analyzing
possible triangle-free graphs arrowing
the corresponding graph parameters
(note that $J_3=P_3$ and thus
any $J_3$-free graph
must be of the form $sK_1 \cup tK_2$).
An example of graph witnessing
$F_v(J_3,J_3;3) \le 9$ can be constructed 
by dropping one vertex from the graph
$C_5[2K_1]$, and for $F_v(K_2,J_3;3) \le 8$
by adding four main diagonals
to the cycle $C_8$. Putting it all together,
by Theorem \ref{FvK{2k+1}-e} applied to this case we obtain
$F_v(J_5,J_5;5) \le 2\cdot 5+2 \cdot 9+8=36$.
\end{proof}

We expect that the actual value of $F_v(J_5,J_5;5)$
is still smaller, but probably not much less so.
How to obtain better bound in this and other similar
cases by detailed analysis is an interesting and
challenging problem.

\medskip
We end this section with two more
upper bounds on $F_v(J_k,J_k;k)$.
Let $E_s$ denote the empty graph on $s$ vertices.
Thus, for example, $K_s[E_t]$ is the same as
the standard complete
$s$-partite graph with all parts of order $t$,
or equivalently, the Tur\'{a}n graph $T_{st,s}$.
One can easily show, similarly as in the proof
of Lemma \ref{prodarrow},
that if $G \rightarrow (k-1,k-1)^v$, then
$G[E_3] \rightarrow (K_{k-1}[E_2],K_{k-1}[E_2])^v$.
Moreover, the same assumption also gives
$G[E_3] \rightarrow (J_k,J_k)^v$,
because $J_k$ is a subgraph of $K_{k-1}[E_2]$.
An even stronger result following from similar
considerations is presented in the next theorem.

\medskip
\begin{theorem} \label{FvKk-e}
If $G \in \mathcal{F}_v(k-1,k-1;k)$,
$|V(G)| = F_v(k-1,k-1;k)$ and $f(G)$
is the largest order of any $K_{k-1}$-free induced
subgraph in $G$, then
$$F_v(J_k,J_k;k) \leq 3 F_v(k-1,k-1;k) - f(G).$$
\end{theorem}

\begin{proof}
Let $G$ be any graph in
$\mathcal{F}_v(k-1, k-1; k)$
with the smallest possible number of vertices
$|V(G)|=F_v(k-1, k-1; k)$. Let us denote
$V(G)=U=\{u_1,\cdots,u_s\}$, so that the vertices
$X=\{u_1,\cdots,u_{f(G)}\}$ induce the largest
$K_{k-1}$-free subgraph in $G$. Note that the
vertices of every $K_{k-1}$ in $G$ have
nonempty intersection with
$Y=\{u_{f(G)+1},\cdots,u_s\}$.
We will construct a graph $H$
on $3|V(G)|-f(G)$ vertices such that
$H \in \mathcal{F}_v(J_k,J_k;k)$, which will
complete the proof of the theorem.

Let $V=\{ v_1,v_2,v_3\}$.
First, take the graph $G[E_3]$ with the set
of vertices $U \times V$. Then $H$ is obtained
from it by dropping $f(G)$ vertices forming the set
$\{\{ u_i,v_3\}\;|\;1 \le i \le f(G)\}$ with
all incident edges. Clearly, graph $H$ has the right
number of vertices and $cl(H)=k-1$. It remains to
be shown that $H \rightarrow (J_k,J_k)^v$.

Let $V(H)=R \cup B$ be any partition
of the vertices of $H$ into two parts,
i.e. any red-blue coloring of $V(H)$.
Note that for each fixed
$u \in U$, there are 2 or 3 vertices in the set
$V(u)=\{(u,v)\;|\;(u,v) \in V(H)\}$.
Let $i(u) \in \{R,B\}$ be a color of at least
half of vertices in $V(u)$ (1 or 2).
Next, consider the red-blue coloring of vertices
of $G$ defined by $i(u)$.
Since $G \rightarrow (k-1,k-1)^v$, then
for the coloring $i(u)$
there exists a set of vertices $S \subset V(G)$
containing a monochromatic $K_{k-1}$ in $G$.
Considering the properties of $X$ and $Y$, $S$ must
contain at least one vertex $u \in Y \cap S$,
and consequently at least two vertices
$(u,x), (u,x) \in V(H)$ are of color $i(u)$.
Now, $S$ expanded to vertices of $H$ of
the same color must contain a monochromatic $J_k$.
\end{proof}

\medskip
\begin{cor}
$$F_v(J_k,J_k;k) \leq
\Bigl\lceil {{5 F_v(k-1, k-1; k)}\over 2}\Bigr\rceil.$$
\end{cor}

\medskip
\begin{proof}
Set $m=\lceil{(F_v(k-1, k-1; k)-1)/2}\rceil$, and let
$G$ and $f(G)$ be as in Theorem \ref{FvKk-e}.
Observe that for every vertex $v \in V(G)$,
$G-v \not\rightarrow (k-1,k-1)^v$ and thus there
exists a $K_ {k-1}$-free set $S \subset V(G) - \{v\}$
satisfying $|S| \ge m$.
Therefore, $f(G) \ge m$ and the corollary easily
follows from Theorem \ref{FvKk-e}.
\end{proof}

\medskip
Using arguments similar to those in
Lemma \ref{prodarrow} and Theorem \ref{FvKk-e},
it is easy to see that
$C_5 [E_{2t-1}] \rightarrow (K_{t,t}, K_{t,t})^v,$
and since $|C_5 [E_{2t-1}]| = 10t-5$, we have an
upper bound
$$F_v (K_{t,t}, K_{t,t}; K_3) \leq 10t -5.$$

\noindent
Improving this bound can be difficult, and
obtaining a good lower bound even harder but
interesting. Hence, we propose the following
problem:

\begin{Prob} \label{F_vK_{t,t}}
For $t \ge 2$,

\smallskip
\noindent
{\rm (a)}
obtain tight bounds for $F_v (K_{t,t}, K_{t,t};3)$, and

\noindent
{\rm (b)}
obtain good bounds for $F_v (K_{t,t}, K_{t,t};k)$, for $k \ge 4$.
\end{Prob}

\medskip
In another application of Theorem \ref{FolkmangeneralGH}
for two colors, with $a_1=a_2=2$, $b=r$ and
$H_1=H_2=T_{tr,r}$,
and for all $t,r \ge 2$, we obtain 
$$F_v(T_{2tr,2r},T_{2tr,2r};2r+1) \leq
5 F_v(T_{tr,r},T_{tr,r}; r+1).$$
In particular,
note that for $r=2$ we have $T_{tr,r}=K_{t,t}$, and
thus the last inequality implies
$F_v(T_{4t,4},T_{4t,4}; K_5) \leq 50t -25.$

\bigskip
The current authors recently studied
the so-called chromatic Folkman numbers \cite{XLR0}, which have
one additional requirement for their witness Folkman
graphs $G$, namely that their chromatic number
$\chi(G)$ is the smallest possible
($\chi(G)=1+ \sum_{i=1}^r (a_i - 1)$
for vertex colorings, and
$\chi(G)=R(a_1, \cdots, a_r)$ for edge colorings).
Some of the constructions in this section lead to
upper bounds involving larger graphs but with
the same chromatic number, mainly because
$\chi(G[E_s]) = \chi(G)$. Thus, these techniques
potentially could lead to stronger claims about upper bounds,
where the chromatic number is as small as possible.

\bigskip
\section{Two Concrete Upper Bounds}

\medskip
If the graphs we wish to arrow to are not of the form
$K_{a_i}[H_i]$ as in Theorem \ref{FolkmangeneralGH},
then the constructions for upper bounds may become
little more complex. For instance, when we deal with
complete but unbalanced multipartite graphs such as
$K_{1,2,2}$ or $K_{2,2,3}$. The cases we study in this section
also involve edge colorings, and the corresponding
generalized edge Folkman numbers. Good upper bounds
for the edge cases seem to be even harder to obtain
than for vertex colorings.
We will focus mainly on two small but puzzling cases
of $F_e(J_4,J_4;5)$ and $F_e(K_3,J_4;5)$.
By the monotonicity of arrowing we have
$$15 = F_e(3,3;5) \le F_e(K_3,J_4;5)
\le F_e(J_4,J_4;5) \le F_e(J_4,J_4;4) \le 30193.$$

\medskip
The equality $F_e(3,3;5)=15$ was obtained
in \cite{Nenov81,PRU}, where in the latter
it was also shown, with the help of computer algorithms,
that $F_v(3,3;4)=14$. Furthermore, in the same work
the authors obtained all 153 graphs on 14 vertices
in the set $\mathcal{F}_v(3,3;4)$. These two parameter
scenarios are connected, since it is known that
$G+u \in \mathcal{F}_e(3,3;5)$ holds whenever
$G \in \mathcal{F}_v(3,3;4)$ (see Lemma \ref{miscarr}(a)).
Two examples of such graphs, $G_a$ and $G_b$,
are presented in Figures A and B.
They are used in the constructions of the following
theorems, which exploit enhancements of
the implication in Lemma \ref{miscarr}(a).
Consult \cite{XLR1}
for the discussion of similar bounds and for
additional pointers to the literature.

\bigskip
\begin{center}
{\small
0 1 1 0 0\ \ 1 1 0 0 0\ \ 1 0 0\ \ 1\\
1 0 0 1 0\ \ 0 0 0 1 1\ \ 0 0 1\ \ 1\\
1 0 0 0 1\ \ 0 0 1 0 1\ \ 0 1 0\ \ 1\\
0 1 0 0 1\ \ 1 0 1 0 0\ \ 0 1 0\ \ 1\\
0 0 1 1 0\ \ 0 1 0 1 0\ \ 0 0 1\ \ 1\\

\smallskip
1 0 0 1 0\ \ 0 0 0 1 1\ \ 1 1 0\ \ 1\\
1 0 0 0 1\ \ 0 0 1 0 1\ \ 1 0 1\ \ 1\\
0 0 1 1 0\ \ 0 1 0 1 0\ \ 1 1 0\ \ 1\\
0 1 0 0 1\ \ 1 0 1 0 0\ \ 1 0 1\ \ 1\\
0 1 1 0 0\ \ 1 1 0 0 0\ \ 0 1 1\ \ 1\\

\smallskip
1 0 0 0 0\ \ 1 1 1 1 0\ \ 0 1 1\ \ 0\\
0 0 1 1 0\ \ 1 0 1 0 1\ \ 1 0 1\ \ 0\\
0 1 0 0 1\ \ 0 1 0 1 1\ \ 1 1 0\ \ 0\\

\smallskip
1 1 1 1 1\ \ 1 1 1 1 1\ \ 0 0 0\ \ 0\\
}
\end{center}

\medskip
\begingroup\leftskip=15pt\rightskip=15pt
\noindent
{\small
{\bf Figure A.} Adjacency matrix of 14-vertex graph
$G_a \in \mathcal{F}_v(3,3;4)$, with a vertex of
maximum degree equal to 10. It is one of 60 such
graphs, all of them enumerated in \cite{PRU}.
Graph $G_a$ has a specially nice structure:
vertices 1--5 and 6--10 induce $C_5$'s,
vertices 11--13 span $K_3$, 10 neighbors of
vertex 14 induce a graph with 320 automorphisms
(which is necessarily triangle-free),
while the entire $G_a$ has only 2 automorphisms.
$G_a$ has independence number 5, it has 41 triangles,
and $|E(G_a)|=48$.
}
\par\endgroup

\eject
\medskip
\begin{center}
{\small
0 0 0 0 0 0 0\ \ \ 1 1 1 1 0 0 0\\ 
0 0 0 0 0 0 0\ \ \ 1 1 1 0 1 0 0\\
0 0 0 0 0 0 0\ \ \ 1 1 0 1 0 1 0\\
0 0 0 0 0 0 0\ \ \ 1 0 1 0 1 0 1\\
0 0 0 0 0 0 0\ \ \ 0 1 0 1 0 1 1\\
0 0 0 0 0 0 0\ \ \ 0 0 1 0 1 1 1\\
0 0 0 0 0 0 0\ \ \ 0 0 0 1 1 1 1\\

\smallskip
1 1 1 1 0 0 0\ \ \ 0 1 1 0 0 1 1\\
1 1 1 0 1 0 0\ \ \ 1 0 0 1 1 0 1\\
1 1 0 1 0 1 0\ \ \ 1 0 0 1 1 1 0\\
1 0 1 0 1 0 1\ \ \ 0 1 1 0 1 1 0\\
0 1 0 1 0 1 1\ \ \ 0 1 1 1 0 0 1\\
0 0 1 0 1 1 1\ \ \ 1 0 1 1 0 0 1\\
0 0 0 1 1 1 1\ \ \ 1 1 0 0 1 1 0\\
}
\end{center}

\medskip
\begingroup\leftskip=15pt\rightskip=15pt
\noindent
{\small
{\bf Figure B.} Adjacency matrix of the Nenov
graph $G_b$ \cite{Nenov81},
which is the unique 14-vertex graph
in the set $\mathcal{F}_v(3,3;4)$
with independence number 7.
In graph $G_b$, vertices 1--7 form
the only 7-independent set and vertices 8--14
induce $\overline{C_7}$.
Graph $G_b$ has 14 automorphisms,
35 triangles, and $|E(G_b)|=42$,
which is the smallest number of edges
among all graphs in $\mathcal{F}_v(3,3;4)$.
}
\par\endgroup

\bigskip
We will also need some simple facts about arrowing.
They are collected in the following lemma.

\smallskip
\begin{lemma} \label{miscarr}
All of the following hold:

\smallskip
\noindent
{\rm (a)}
if $G\rightarrow (3,3)^v$ and $u$ is a new vertex,
then $G+u \rightarrow (3,3)^e$,

\smallskip
\noindent
{\rm (b)}
if $G\rightarrow (3,3)^v$, then
$G[E_{2k-1}] \rightarrow (K_{k,k,k},K_{k,k,k})^v$ for $k \ge 1$,

\smallskip
\noindent
{\rm (c)}
$K_{1,2,2} \rightarrow (J_3,K_3)^e$, and

\smallskip
\noindent
{\rm (d)}
$K_{2,2,3} \rightarrow (J_3,J_4)^e$.
\end{lemma}

\begin{proof}
Part (a) is a basic property of arrowing
used by many authors in scenarios similar to ours,
for example, in \cite{Nenov81}.
For part (b) observe that $K_3[E_k]=K_{k,k,k}$
and use Lemma \ref{prodarrow} for two colors with
$r=2$, $a_1=a_2=3$ and $H_1=H_2=E_k$. For parts
(c) and (d), first note that
$K_{1,2,2}=K_1+C_4$ and $K_{2,2,3}=E_3+C_4$,
and then consider possible edges of $C_4$ in the
color avoiding $J_3$; there are at most two such edges.
There are just a few more choices of edges in the
$J_3$-free color since in total there are at most
2 or 3 such edges in $K_{1,2,2}$ and $K_{2,2,3}$,
respectively. A routine check easily shows that
the edges in the other color must contain a $K_3$
and $J_4$, respectively.
\end{proof}

\medskip
\begin{theorem} \label{FvJ4J4K5}
$F_e(J_4,J_4;5) \leq 51$.
\end{theorem}

\begin{proof}
The skeleton of the proof is as follows.
First, we will construct a $K_4$-free graph $H$
on 50 vertices such that
$H \rightarrow (K_{2,2,3},K_{2,2,3})^v$.
Then we will claim that the 51-vertex graph
$G=K_1+H$ is in the set $\mathcal{F}_e(J_4,J_4;5)$,
or equivalently, it is a witness of the upper bound
in the theorem.

We start with the graph $G_a \in \mathcal{F}_v(3,3;4)$
described in Figure A.
Let $u$ be the only vertex of degree 10.
Consider the partition of the set $V(G_a)$ into
$X \cup Y$, where
$X=\{v\;|\;\{u,v\}\in E(G_a)\}$,
so that $|X|=10$ and $|Y|=4$.
Note that every triangle in $G_a$ must have
at least one vertex in $Y$.

Next, define the graph $G_1=G_a[E_5]$ which
has 70 vertices and it is $K_4$-free.
By Lemma \ref{miscarr}(b), it holds that
$G_1 \rightarrow (K_{3,3,3},K_{3,3,3})^v$,
which is too strong for our purpose, so we can
drop some vertices from $G_1$.
We define graph $H$ as an induced subgraph
of $G_1$ following the idea of the proof of
Theorem \ref{FvKk-e}, which now is to
drop $2|X|$ vertices from $G_1$ formed
by two triangle-free parts of $G_a$.
More precisely, if the product graph $G_1$ uses
$V(E_5)=\{v_1,\cdots,v_5\}$, then
$V(H)=V(G_1) \setminus X \times \{v_4,v_5\}$
and $|V(H)|=50$.
By an argument very similar to that in the proof
of Theorem \ref{FvKk-e} we can see that
$H \rightarrow (K_{2,2,3},K_{2,2,3})^v$.
Being a subgraph of $G_1$, the graph $H$ is $K_4$-free.

Finally, define graph $G$ to be $K_1+H$, and
let $u$ be the vertex in $V(G) \setminus V(H)$.
Clearly, $G$ is $K_5$-free and we have
$V(G)=51$. It remains to be shown that
$G \rightarrow (J_4,J_4)^e$. Consider any
red-blue coloring $C_e$ of the edges $E(G)$.
Define a red-blue vertex coloring $C_v$
by $C_v(x)=C_e(\{u,x\})$ for $x \in V(H)$.
Since $H \rightarrow (K_{2,2,3},K_{2,2,3})^v$,
then $H$ contains a monochromatic $K_{2,2,3}$
in $C_v$, say, with the vertex set $U$.
Without loss of generality assume that
all vertices in $U$ are red.
Now, by Lemma \ref{miscarr}(d),
the set of edges induced by $U$ must contain
a red $J_3$ or blue $J_4$ in $C_e$.
If it is blue $J_4$, then we are done, otherwise
red $J_3$ together with vertex
$u$ induces a red $J_4$.
\end{proof}

\medskip
\begin{theorem} \label{FvK3J4K5}
$F_e(K_3,J_4; K_5) \leq 27$.
\end{theorem}

\begin{proof}
The reasoning is very similar to that in the proof of
Theorem \ref{FvJ4J4K5}, just the parameters vary.
We start with the graph $G_b$ described in Figure B,
and let $G_1=G_b[E_3]$.
By Lemma \ref{miscarr}(b), it holds that
$G_1 \rightarrow (K_{2,2,2},K_{2,2,2})^v$,
so we can drop some vertices from $G_1$.
Note that $J_4=K_{1,1,2}$.
We will define an induced subgraph $H$ of $G_1$
on 26 vertices
such that $H \rightarrow (K_{1,1,2},K_{1,2,2})^v$.
Then the graph $K_1+H$ will be a witness
of the upper bound.

Consider a partition of the set $V(G_b)$ into
$X \cup Y \cup Z$, where
$X$ consists of the vertices of 7-independent set
(the first 7 vertices), $Y$ is the pair of vertices
of any nonedge contained in the second block of 7 vertices,
$Z$ is formed by the remaining vertices,
and thus they have orders 7, 2 and 5, respectively.
Note that every triangle in $G_b$ must have
at least one vertex in $Z$, and at least
two vertices in $Y \cup Z$.
Drop $2|X|$ vertices from $V(G_1)$, 2 associated with
each vertex in $X$, and similarly drop $|Y|$ vertices
associated with $Y$. This defines an induced
subgraph $H$ on $|X|+2|Y|+3|Z|=26$ vertices.

Now,
for every red-blue vertex coloring $C_1$ of $V(G_1)$
define the coloring $C_b$ of $V(G_b)$ by
assigning to $u \in V(G_b)$ color
of the majority of vertices in the set
$\bigl\{\{u,v\}\;|\;\{u,v\}\in V(H)\bigr\}$ under $C_1$.
Note that these sets have the cardinality 1, 2 and 3
for $u$ in $X$, $Y$ and $Z$, respectively.
Thus, every monochromatic triangle in $G_b$ can be
expanded in $H$ to at least 1, 1, and 3 vertices, or
at least 1, 2, and 2 vertices, respectively.
Hence, similarly as in the proof of
Theorem \ref{FvKk-e}, we can see that
$H \rightarrow (K_{1,1,2},K_{1,2,2})^v$. Furthermore,
being a subgraph of $G_1$, the graph $H$ is $K_4$-free.

Finally, define the graph $G$ as $K_1+H$,
so that $G$ is $K_5$-free and we have
$V(G)=27$. It remains to be shown that
$G \rightarrow (K_3,J_4)^e$.
Let $u$ be the vertex in $V(G) \setminus V(H)$
and consider any
red-blue coloring $C_e$ of the edges $E(G)$.
Define a red-blue vertex coloring $C_v$
by $C_v(x)=C_e(\{u,x\})$ for $x \in V(H)$.
Since $H \rightarrow (K_{1,1,2},K_{1,2,2})^v$,
then $H$ contains a red $K_{1,1,2}$ or blue
$K_{1,2,2}$ in $C_v$, say with the vertex set $U$.
Recall that $J_4=K_{1,1,2}$. 
First suppose that
the vertices in $U$ are forming a red $J_4$ in $C_v$.
If at least one of the edges with both endpoints
in $U$ is red in $C_e$, then we have red $K_3$
spanned by this edge and $u$ in $C_e$. Otherwise,
all edges induced by $U$ are blue in $C_e$, so we
have blue $J_4$. On the other hand, suppose that
the vertices in $U$ form a blue $K_{1,2,2}$ in $C_v$.
By Lemma \ref{miscarr}(c),
the set of edges induced by $U$ must contain
a red $K_3$ or blue $J_3$ in $C_e$.
If it is red $K_3$, then we are done, otherwise
blue $J_3$ together with vertex
$u$ induces a blue $J_4$.
\end{proof}

\medskip
The upper bounds of the last two theorems are
likely larger than the actual values, despite that
the graphs $G_a$ and $G_b$ where chosen
to be the best for our method, namely: none of the graphs on
14 vertices in $\mathcal{F}_v(3,3;4)$ has any induced triangle-free
subgraph on more than 10 vertices and none has two independent
sets whose union has more than 9 vertices, respectively.
On the other hand, improving on either of the bounds
$F_v(K_{2,2,3},K_{2,2,3};4) \le 50$ or
$F_v(K_{1,1,2},K_{1,2,2};4) \le 26$, currently
used in the proofs, would lead to better upper bounds
in Theorems 10 and 11, respectively. Finally, let us
note that Folkman numbers with other parameters
could be studied by the current method, for example,
such as when exploiting the inequality
$F_e(3,4;7) \leq 1+F_v(K_4,K_3+C_5;6)$, which relies on
the well known case of edge arrowing $K_3+C_5 \rightarrow (3,3)^e$.

\end{document}